# A Class of Orthohomological Triangles

Prof. Claudiu Coandă, National College "Carol I", Craiova, Romania;
E-mail: coandac@yahoo.com.
Prof. Florentin Smarandache, University of New Mexico, Gallup, USA;
E-mail: smarand@unm.edu.
Prof. Ion Pătrașcu, National College "Frații Buzești", Craiova, Romania.
E-mail: patrascu_ion@yahoo.com.

**Abstract.**
In this article we propose to determine the triangles' class $A_iB_iC_i$ orthohomological with a given triangle $ABC$, inscribed in the triangle $ABC$ ($A_i \in BC$, $B_i \in AC$, $C_i \in AB$).

We'll remind, here, the fact that if the triangle $A_iB_iC_i$ inscribed in $ABC$ is orthohomologic with it, then the perpendiculars in $A_i$, $B_i$, respectively in $C_i$ on $BC$, $CA$, respectively $AB$ are concurrent in a point $P_i$ (the orthological center of the given triangles), and the lines $AA_i$, $BB_i$, $CC_i$ are concurrent in point (the homological center of the given triangles).

To find the triangles $A_iB_iC_i$, it will be sufficient to solve the following problem.

**Problem.**
Let's consider a point $P_i$ in the plane of the triangle $ABC$ and $A_iB_iC_i$ its pedal triangle. Determine the locus of point $P_i$ such that the triangles $ABC$ and $A_iB_iC_i$ to be homological.

**Solution.**
Let's consider the triangle $ABC$, $A(1,0,0)$, $B(0,1,0)$, $C(0,0,1)$, and the point $P_i(\alpha, \beta, \gamma)$, $\alpha + \beta + \gamma = 0$.

The perpendicular vectors on the sides are:
$$U^\perp_{BC}\left(2a^2,\ -a^2-b^2+c^2,\ -a^2+b^2-c^2\right)$$
$$U^\perp_{CA}\left(-a^2-b^2+c^2,\ 2b^2,\ a^2-b^2-c^2\right)$$
$$U^\perp_{AB}\left(-a^2+b^2-c^2,\ a^2-b^2-c^2,\ 2c^2\right)$$

The coordinates of the vector $\overrightarrow{BC}$ are $(0,-1,1)$, and the line $BC$ has the equation $x=0$.

The equation of the perpendicular raised from point $P_i$ on $BC$ is:

$$\begin{vmatrix} x & y & z \\ \alpha & \beta & \gamma \\ 2a^2 & -a^2-b^2+c^2 & -a^2+b^2-c^2 \end{vmatrix} = 0$$

We note $A_i(x,y,z)$, because $A_i \in BC$ we have:

$x = 0$ and $y + z = 1$.

The coordinates $y$ and $z$ of $A_i$ can be found by solving the system of equations

$$\begin{cases} \begin{vmatrix} x & y & z \\ \alpha & \beta & \gamma \\ 2a^2 & -a^2-b^2+c^2 & -a^2+b^2-c^2 \end{vmatrix} = 0 \\ y+z = 0 \end{cases}$$

We have:

$$y \cdot \begin{vmatrix} \alpha & \gamma \\ 2a^2 & -a^2+b^2-c^2 \end{vmatrix} = z \cdot \begin{vmatrix} \alpha & \beta \\ 2a^2 & -a^2-b^2+c^2 \end{vmatrix},$$

$$y\left[\alpha\left(-a^2+b^2-c^2\right)-2\gamma a^2\right] = z\left[\alpha\left(-a^2-b^2+c^2\right)-2\beta a^2\right],$$

$$y + y \cdot \frac{\alpha\left(-a^2+b^2-c^2\right)-2\gamma a^2}{\alpha\left(-a^2-b^2+c^2\right)-2\beta a^2} = 1,$$

$$y \cdot \frac{\alpha\left(-a^2-b^2+c^2\right)-2\beta a^2+\alpha\left(-a^2+b^2-c^2\right)-2\gamma a^2}{\alpha\left(-a^2-b^2+c^2\right)-2\beta a^2} = 1,$$

$$y \cdot \frac{-2a^2(\alpha+\beta+\gamma)}{\alpha\left(-a^2-b^2+c^2\right)-2\beta a^2} = 1,$$

it results

$$y = \frac{\alpha}{2a^2}\left(a^2+b^2-c^2\right) + \beta$$

$$z = 1 - y = 1 - \beta - \frac{\alpha}{2a^2}\left(a^2+b^2-c^2\right) = \alpha + \gamma - \frac{\alpha}{2a^2}\left(a^2+b^2-c^2\right).$$

Therefore,

$$A_i\left(0,\ \frac{\alpha}{2a^2}\left(a^2+b^2-c^2\right)+\beta,\ \frac{\alpha}{2a^2}\left(a^2-b^2+c^2\right)+\gamma\right).$$

Similarly we find:

$$B_i\left(\frac{\beta}{2b^2}\left(a^2+b^2-c^2\right)+\alpha,\ 0,\ \frac{\beta}{2b^2}\left(-a^2+b^2+c^2\right)+\gamma\right),$$

$$C_i\left(\frac{\gamma}{2c^2}\left(a^2-b^2+c^2\right)+\alpha,\ \frac{\gamma}{2c^2}\left(-a^2+b^2+c^2\right)+\beta,\ 0\right).$$

We have:

$$\overline{\dfrac{A_iB}{A_iC}} = -\dfrac{\dfrac{\alpha}{2a^2}(a^2-b^2+c^2)+\gamma}{\dfrac{\alpha}{2a^2}(a^2+b^2-c^2)+\beta} = -\dfrac{\alpha c\cos B+\gamma a}{\alpha b\cos C+\beta a}$$

$$\overline{\dfrac{B_iC}{B_iA}} = -\dfrac{\dfrac{\beta}{2b^2}(a^2+b^2-c^2)+\alpha}{\dfrac{\alpha}{2a^2}(-a^2+b^2+c^2)+\gamma} = -\dfrac{\beta a\cos C+\alpha b}{\beta c\cos A+\gamma b}.$$

$$\overline{\dfrac{C_iA}{C_iB}} = -\dfrac{\dfrac{\gamma}{2c^2}(-a^2+b^2+c^2)+\beta}{\dfrac{\gamma}{2c^2}(a^2-b^2+c^2)+\alpha} = -\dfrac{\gamma b\cos A+\beta c}{\gamma a\cos B+\alpha c}$$

(We took into consideration the cosine's theorem: $a^2 = b^2+c^2-2bc\cos A$).
In conformity with Ceva's theorem, we have:

$$\overline{\dfrac{A_iB}{A_iC}} \cdot \overline{\dfrac{B_iC}{B_iA}} \cdot \overline{\dfrac{C_iA}{C_iB}} = -1.$$

$$(a\gamma+\alpha c\cos B)(b\alpha+\beta a\cos C)(c\beta+\gamma b\cos A) =$$
$$= (a\beta+\alpha b\cos C)(b\gamma+\beta c\cos A)(c\alpha+\gamma a\cos B)$$
$$a\alpha(b^2\gamma^2-c^2\beta^2)(\cos A-\cos B\cos C)+b\beta(c^2\alpha^2-a^2\gamma^2)(\cos B-\cos A\cos C)+$$
$$+c\gamma(a^2\beta^2-b^2\alpha^2)(\cos C-\cos A\cos B) = 0.$$

Dividing it by $a^2b^2c^2$, we obtain that the equation in barycentric coordinates of the locus $\mathcal{L}$ of the point $P_i$ is:

$$\dfrac{\alpha}{a}\left(\dfrac{\gamma^2}{c^2}-\dfrac{\beta^2}{b^2}\right)(\cos A-\cos B\cos C)+\dfrac{\beta}{b}\left(\dfrac{\alpha^2}{a^2}-\dfrac{\gamma^2}{c^2}\right)(\cos B-\cos A\cos C)+$$
$$+\dfrac{\gamma}{c}\left(\dfrac{\beta^2}{b^2}-\dfrac{\alpha^2}{a^2}\right)(\cos C-\cos A\cos B) = 0.$$

We note $\bar{d}_A$, $\bar{d}_B$, $\bar{d}_C$ the distances oriented from the point $P_i$ to the sides $BC$, $CA$ respectively $AB$, and we have:

$$\dfrac{\alpha}{a} = \dfrac{\bar{d}_A}{2s},\ \dfrac{\beta}{b} = \dfrac{\bar{d}_B}{2s},\ \dfrac{\gamma}{c} = \dfrac{\bar{d}_C}{2s}.$$

The locus' $\mathcal{L}$ equation can be written as follows:

$$\bar{d}_A(\bar{d}_C^2-\bar{d}_B^2)(\cos A-\cos B\cos C)+\bar{d}_B(\bar{d}_A^2-\bar{d}_C^2)(\cos B-\cos A\cos C)+$$
$$+\bar{d}_C(\bar{d}_B^2-\bar{d}_A^2)(\cos C-\cos A\cos B) = 0$$

**Remarks.**

1. It is obvious that the triangle's $ABC$ orthocenter belongs to locus $\mathcal{L}$. The orthic triangle and the triangle $ABC$ are orthohomologic; a orthological center is the orthocenter $H$, which is the center of homology.
2. The center of the inscribed circle in the triangle $ABC$ belongs to the locus $\mathcal{L}$, because $\bar{d}_A = \bar{d}_B = \bar{d}_C = r$ and thus locus' equation is quickly verified.

**Theorem** *(Smarandache-Pătraşcu).*
If a point $P$ belongs to locus $\mathcal{L}$, then also its isogonal $P'$ belongs to locus $\mathcal{L}$.

**Proof.**
Let $P(\alpha, \beta, \gamma)$ a point that verifies the locus' $\mathcal{L}$ equation, and $P'(\alpha', \beta', \gamma')$ its isogonal in the triangle $ABC$. It is known that $\dfrac{\alpha \alpha'}{a^2} = \dfrac{\beta \beta'}{b^2} = \dfrac{\gamma \gamma'}{c^2}$. We'll prove that $P' \in \mathcal{L}$, i.e.

$$\sum \frac{\alpha'}{a}\left(\frac{\gamma'^2}{c^2} - \frac{\beta'^2}{b^2}\right)(\cos A - \cos B \cos C) = 0$$

$$\sum \frac{\alpha'}{a}\left(\frac{\gamma'^2 b^2 - \beta'^2 c^2}{b^2 c^2}\right)(\cos A - \cos B \cos C) = 0$$

$$\sum \frac{\alpha'}{ab^2 c^2}\left(\gamma'^2 b^2 - \beta'^2 c^2\right)(\cos A - \cos B \cos C) = 0 \Leftrightarrow$$

$$\Leftrightarrow \sum \frac{\alpha'}{ab^2 c^2}\left(\frac{\gamma' \beta \beta' c^2}{\gamma} - \frac{c^2 \gamma \gamma' \beta'}{\beta}\right)(\cos A - \cos B \cos C) = 0 \Leftrightarrow$$

$$\Leftrightarrow \sum \frac{\alpha' \beta' \gamma'}{ab^2 c^2}\left(\frac{\beta c^2}{\gamma} - \frac{\gamma b^2}{\beta}\right)(\cos A - \cos B \cos C) = 0 \Leftrightarrow$$

$$\Leftrightarrow \sum \frac{\alpha' \beta' \gamma'}{ab^2 c^2}\left(\frac{\beta^2 c^2 - \gamma^2 b^2}{\beta \gamma}\right)(\cos A - \cos B \cos C) = 0 \Leftrightarrow$$

$$\Leftrightarrow \sum \frac{\alpha}{a}\left(\frac{\alpha' \beta' \gamma'}{\alpha \beta \gamma}\right)\frac{1}{b^2 c^2} \cdot b^2 c^2 \left(\frac{\beta^2}{b^2} - \frac{\gamma^2}{c^2}\right)(\cos A - \cos B \cos C) = 0.$$

We obtain that:

$$\frac{\alpha' \beta' \gamma'}{\alpha \beta \gamma} \sum \frac{\alpha}{a}\left(\frac{\gamma^2}{c^2} - \frac{\beta^2}{b^2}\right)(\cos A - \cos B \cos C) = 0,$$

this is true because $P \in \mathcal{L}$.

**Remark.**
We saw that the triangle's $ABC$ orthocenter $H$ belongs to the locus, from the precedent theorem it results that also $O$, the center of the circumscribed circle to the triangle $ABC$ (isogonable to $H$), belongs to the locus.

**Open problem:**
What does it represent from the geometry's point of view the equation of locus $\mathcal{L}$?

In the particular case of an equilateral triangle we can formulate the following:

**Proposition:**
The locus of the point $P$ from the plane of the equilateral triangle $ABC$ with the property that the pedal triangle of $P$ and the triangle $ABC$ are homological, is the union of the triangle's heights.

**Proof:**
Let $P(\alpha, \beta, \gamma)$ a point that belongs to locus $\mathcal{L}$. The equation of the locus becomes:
$$\alpha(\gamma^2 - \beta^2) + \beta(\alpha^2 - \gamma^2) + \gamma(\beta^2 - \alpha^2) = 0$$
Because:
$$\alpha(\gamma^2 - \beta^2) + \beta(\alpha^2 - \gamma^2) + \gamma(\beta^2 - \alpha^2) = \alpha\gamma^2 - \alpha\beta^2 + \beta\alpha^2 - \beta\gamma^2 + \gamma\beta^2 - \gamma\alpha^2 =$$
$$= \alpha\beta\gamma + \alpha\gamma^2 - \alpha\beta^2 + \beta\alpha^2 - \beta\gamma^2 + \gamma\beta^2 - \gamma\alpha^2 - \alpha\beta\gamma =$$
$$= \alpha\beta(\gamma - \beta) + \alpha\gamma(\gamma - \beta) - \alpha^2(\gamma - \beta) - \beta\gamma(\gamma - \beta) =$$
$$= (\gamma - \beta)\left[\alpha(\beta - \alpha) - \gamma(\beta - \alpha)\right] = (\beta - \alpha)(\alpha - \gamma)(\gamma - \beta).$$

We obtain that $\alpha = \beta$ or $\beta = \gamma$ or $\gamma = \alpha$, that shows that $P$ belongs to the medians (heights) of the triangle $ABC$.